\documentclass[12pt]{article}
\usepackage{latexsym,amsmath, amsfonts, amscd, amsthm, mathrsfs}
\usepackage{amssymb}
\usepackage{mathptm}    
\usepackage[all,dvips,arc,curve,color,frame]{xy}

\newtheorem{theorem}{Theorem}[section]
\newtheorem{lemma}[theorem]{Lemma}
\newtheorem{proposition}[theorem]{Proposition}
\newtheorem{corollary}[theorem]{Corollary}
\def\square{\Box}

\newenvironment{prf}[1]{\trivlist
\item[\hskip
\labelsep{\it #1.\hspace*{.3em}}]}{~\hspace{\fill}~$\square$
\endtrivlist}

\newtheorem{predefinition}[theorem]{Definition}

\newtheorem{preremark}[theorem]{Remark}
\newenvironment{remark}{\begin{preremark}\rm}{\end{preremark}}
\newtheorem{prenotation}[theorem]{Notation}
\newenvironment{notation}{\begin{prenotation}\rm}{\end{prenotation}}
\newtheorem{preexample}[theorem]{Example}
\newenvironment{example}{\begin{preexample}\rm}{\end{preexample}}
\newtheorem{preclaim}[theorem]{Claim}

\newtheorem{prequestion}[theorem]{Question}

\theoremstyle{remark}
\newtheorem{acknowledgments}{Acknowledgments}

\newcommand \CA {{\mathcal A}}
\newcommand \CM {{\mathcal M}}
\newcommand{\AS}{{\mathcal{AS}}}
\newcommand {\ArSc} {{\mathcal {AS}cov}}
\newcommand \CH {{\mathcal H}}
\newcommand \CO {{\mathcal O}}
\newcommand \Cdiv {{\mathcal Cdiv}}

\newcommand \PP {{\mathbb P}^1}
\newcommand \ZZ {{\mathbb Z}}

\newcommand \Aa {{\mathbb A}^1}

\newcommand \Aut {\mathop{\rm Aut}}
\newcommand \dime {\mathop{\rm dim}}
\newcommand \Jac {\mathop{\rm Jac}}
\newcommand \Spec {\mathop{\rm Spec}}

\title{The $p$-rank stratification of Artin-Schreier curves}

\author{Rachel Pries and Hui June Zhu
\footnote{The first author was partially supported by NSF grant DMS 07-01303.
The second author was partially supported by an NSA grant.}}

\begin{document}

\maketitle

\begin{abstract}
We study a moduli space  ${\mathcal{AS}}_g$ for Artin-Schreier curves of genus $g$
over an algebraically closed field $k$ of characteristic $p$.
We study the stratification of ${\mathcal{AS}}_g$ by $p$-rank into strata ${\mathcal{AS}}_{g.s}$
of Artin-Schreier curves of genus $g$ with $p$-rank exactly $s$.
We enumerate the irreducible components of ${\mathcal{AS}}_{g,s}$ and find their dimensions.
As an application, when $p=2$, we prove that every irreducible component of the
moduli space of hyperelliptic $k$-curves with genus $g$ and $2$-rank $s$ has dimension $g-1+s$.
We also determine all pairs $(p,g)$ for which ${\mathcal{AS}}_g$ is irreducible.
Finally, we study deformations of Artin-Schreier curves with varying $p$-rank.\\
Keywords: Artin-Schreier, hyperelliptic, curve, moduli, $p$-rank.\\
MSC: 11G15, 14H40, 14K15.
\end{abstract}

\section{Introduction}

Let $k$ be an algebraically closed field of characteristic $p >0$.
An {\it Artin-Schreier $k$-curve} is a smooth projective
connected $k$-curve $Y$ which is a $(\ZZ/p)$-cover of the
projective line. The Riemann-Hurwitz formula implies that the genus $g$ of $Y$ is of the
form $g=d(p-1)/2$ for some integer $d \geq 0$. The {\it
$p$-rank} of $Y$ is the integer $s$ such
that the cardinality of $\Jac(Y)[p](k)$ is $p^s$.
It is well known that $0 \leq s \leq g$.
By the Deuring-Shafarevich formula, $s=r(p-1)$ for some integer $r \geq 0$.

In this paper, we study a moduli space $\AS_{g}$ for Artin-Schreier $k$-curves of genus $g$.
We study its stratification by $p$-rank into strata $\AS_{g,s}$ whose points
correspond to Artin-Schreier curves of genus $g$ with $p$-rank exactly $s$.
Throughout, we assume $g=d(p-1)/2$ and $s=r(p-1)$ for some integers $d \geq 1$ and $r \geq 0$ since
the problem is trivial otherwise.
We denote by $\lfloor\cdot \rfloor$ and $\lceil\cdot\rceil$ the floor and
ceiling of a real number, respectively, and use the notation $\{\cdots\}$ to denote a multi-set.
We prove:

\begin{theorem} \label{T1}
Let $g=d(p-1)/2$ with $d \geq 1$ and $s=r(p-1)$ with $r \geq 0$.
\begin{enumerate}
\item The set of irreducible components of $\AS_{g,s}$ is in
bijection with the set of partitions $\{e_1, \ldots e_{r+1}\}$ of
$d+2$ into $r+1$ positive integers such that each $e_j \not \equiv 1 \bmod p$.
\item The irreducible component of $\AS_{g,s}$
for the partition $\{e_1, \ldots e_{r+1}\}$ has dimension
$$d-1 -\sum_{j=1}^{r+1} \lfloor (e_j-1)/p \rfloor.$$
\end{enumerate}
\end{theorem}

The proof uses ideas from \cite[Section 5.1]{BM}, \cite{Ha:mod}, and \cite{Pr:fam}.
As an application of Theorem \ref{T1}, we determine all cases when $\AS_g$ is irreducible,
using the fact that every irreducible component of $\AS_g$ has dimension $d-1$, \cite[Cor.\ 3.16]{Maug}.

\begin{corollary} \label{C2}
The moduli space $\AS_g$ is irreducible in exactly the following cases:
(i) $p=2$; or (ii) $g=0$ or $g=(p-1)/2$; or (iii) $p=3$ and $g=2,3,5$.
\end{corollary}

When $p=2$, the moduli space $\AS_g$ is the same as $\CH_g$, the moduli space of hyperelliptic $k$-curves
of genus $g$.  By \cite[Thm.\ 4.1]{Lonsted},
$\CH_g$ is irreducible of dimension $2g-1$ when $p=2$.
Let $\CH_{g,s} \subset \CH_g$ denote the stratum whose points correspond to
hyperelliptic $k$-curves of genus $g$ with $2$-rank $s$.
Theorem \ref{T1} yields the following description of $\CH_{g,s}$.
This also generalizes the result $\dime(\CH_{g,0})=g-1$ when $p=2$ from \cite[Prop.\ 4.1]{SZ:02}.

\begin{corollary} \label{C1}
Let $p=2$ and $g \geq 1$.  The irreducible components of $\CH_{g,s}$ are in
bijection with partitions of $g+1$ into $s+1$ positive integers.
Every component has dimension $g-1+s$.
\end{corollary}

The geometry of $\AS_g$ is more complicated when $p \geq 3$.
For example, Theorem \ref{T1} shows that, for fixed $g$ and $s$, the irreducible components of $\AS_{g,s}$
can have different dimensions and thus $\AS_{g,s}$ is not pure in general when $p \geq 3$, Corollary \ref{Cnottrans}.

Here is some motivation for these results, which also gives another illustration how the geometry of
$\AS_g$ is more complicated when $p \geq 3$.
Recall that the moduli space $\CA_g$
of principally polarized abelian varieties over $k$ of dimension $g$ can be
stratified by $p$-rank. Let $V_{g,s} \subset \CA_g$ denote the stratum of
abelian varieties with $p$-rank $s$. By \cite[1.6]{O:purity}, every component of $V_{g,s}$
has codimension $g-s$ in $\CA_g$.
Suppose $M$ is a subspace of $\CM_g$, the moduli space of $k$-curves of genus $g$.
One can ask whether the image $T(M)$ of $M$ under the Torelli morphism is in general position relative to the $p$-rank
stratification. A necessary condition for an affirmative answer is
that ${\rm codim}(T(M) \cap V_{g,s}, T(M))=g-s$. This has been verified
when $M=\CM_g$ in \cite[Thm.\ 2.3]{FVdG:complete}
and when $M=\CH_g$ for $p \geq 3$ in \cite[Thm.\ 1]{GP:05}.
Corollary \ref{C1} shows that this necessary condition is satisfied for $M=\CH_g$ when $p=2$.
Corollary \ref{Cnottrans} shows that it is not satisfied for $M=\AS_g$ when $p\geq 3$.

Finally, we study how the components of $\AS_{g,s}$ (with varying $s$) fit together inside $\AS_g$.
This is related to the study of deformations of wildly ramified degree $p$ covers with non-constant branch locus. Under the obvious necessary conditions, we prove that the $p$-rank of an Artin-Schreier curve can
be increased by exactly $p-1$ in a flat deformation.  This yields the following result.

\begin{theorem} \label{Tclosure}
Suppose $0 \leq s \leq g-(p-1)$.
If $\eta$ is an irreducible component of $\AS_{g,s}$ which is not
open and dense in an irreducible component of $\AS_g$,
then $\eta$ is in the closure of $\AS_{g,s+(p-1)}$ in $\AS_g$.
\end{theorem}

When $p=2$, we are further able to give
a complete combinatorial description of how the irreducible components of
$\CH_{g,s}$ (with varying $s$) fit together in $\CH_g$, Corollary \ref{deform2}.

Here is an outline of the paper. In Section \ref{S2}, we describe
the $p$-ranks of Artin-Schreier curves and the relationship between irreducible components and partitions.
Section \ref{S4} contains the proof of the main results.
One finds Theorem \ref{T1} in
Section \ref{Scover}, Corollary \ref{C2} in Section \ref{Sirred},
and Corollary \ref{C1} in Section \ref{Shyper}.
The deformation results, including Theorem \ref{Tclosure}, are in Section \ref{Sdeform}.
We conclude with some open questions.

\begin{acknowledgments}
The authors thank the referee, Jeff Achter, David Harbater, and Jasper Scholten
for comments and correspondences about this paper.
The second author thanks Rachel Pries and Jeff Achter for
invitation and hospitality during her visit to Colorado, where this work was initiated.
\end{acknowledgments}

\section{Partitions and Artin-Schreier curves} \label{S2}

\subsection{Partitions} \label{S3}

Fix a prime $p>0$ and an integer $d\geq 1$, with $d$ even if $p=2$.  Let $\Omega_{d}$ be the set of
partitions of $d+2$ into positive integers $e_1,e_2,\ldots$ with
each $e_j \not \equiv 1 \bmod p$. Let $\Omega_{d,r}$ be the subset of
$\Omega_d$ consisting of partitions of length $r+1$.
If $\vec{E} \in \Omega_d$, let $r:=r(\vec{E})$ be the integer so
that $\vec{E} \in \Omega_{d,r}$.
Write $\vec{E}$ as a multi-set $\{e_1, \ldots, e_{r+1}\}$ with $e_1 \leq \cdots \leq e_{r+1}$.

There is a natural partial ordering $\prec$ on $\Omega_d$ so that
$\vec{E} \prec \vec{E}'$ if $\vec{E}'$ is a refinement of
$\vec{E}$, in other words, if the entries of $\vec{E}'$ can be
divided into disjoint subsets whose sums are in bijection with the
entries of $\vec{E}$. Using this partial ordering, one can
construct a directed graph $G_d$. The vertices of the graph
correspond to the partitions $\vec{E}$ in $\Omega_d$. There is an
edge from $\vec{E}$ to $\vec{E}'$ if and only if $\vec{E} \prec
\vec{E}'$, and $\vec{E} \not = \vec{E}'$, and there is no
partition lying strictly in between them (i.e., if $\vec{E} \prec
\vec{E''} \prec \vec{E'}$ for some $\vec{E''}$ in $\Omega_d$ then
$\vec{E''} = \vec{E}$ or $\vec{E''} = \vec{E'}$).

An edge $\vec{E} \prec \vec{E'}$ in the directed graph $G_d$ can
be of two types.  The first type has $r(\vec{E})= r(\vec{E}')-1$.
In this case, one entry $e$ of $\vec{E}$ splits into two entries
$e_1$ and $e_2$ of $\vec{E}'$ such that $e = e_1 +e_2$ and none of
the three is congruent to $1$ modulo $p$. One can summarize this
by writing $\{e\} \mapsto \{e_1,e_2\}$. The second type has
$r(\vec{E})= r(\vec{E}')-2$. In this case, one entry $e$ of
$\vec{E}$ splits into three entries $e_1$, $e_2$, $e_3$ of
$\vec{E}'$ such that $e=e_1 +e_2+e_3$ and each $e_j \equiv (p+1)/2
\bmod p$. It follows that none of the four is congruent to $1$
modulo $p$. One can summarize this by writing $\{e\} \mapsto
\{e_1,e_2,e_3\}$.

\begin{example}\label{E:1}
Let $p=3$ and $d=10$.
Here is the graph $G_{10}$ for $\Omega_{10}$.
$$
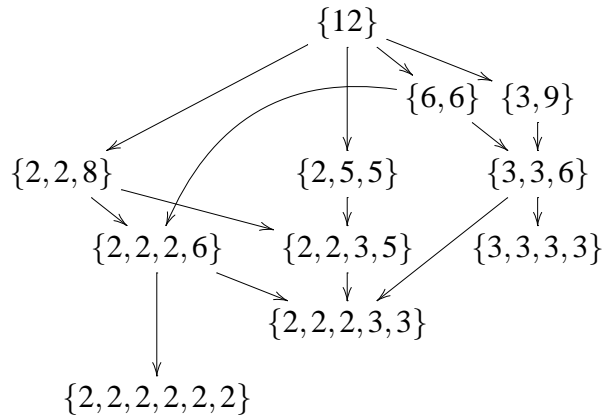
\xygraph{
!{<0cm,0cm>;<1cm,0cm>;<0cm,1cm>::}
!{(0,5)*+{{\{12\}}}}="a"
!{(2,4)*+{{\{3,9\}}}}="b1"
!{(1,4)*+{{\{6,6\}}}}="b2"
!{(-3,3)*+{{\{2,2,8\}}}}="c1"
!{(0,3)*+{{\{2,5,5\}}}}="c2"
!{(2,3)*+{{\{3,3,6\}}}}="c3"
!{(-2,2)*+{{\{2,2,2,6\}}}}="d1"
!{(0,2)*+{{\{2,2,3,5\}}}}="d2"
!{(2,2)*+{{\{3,3,3,3\}}}}="d3"
!{(0,1)*+{{\{2,2,2,3,3\}}}}="e"
!{(-2,0)*+{{\{2,2,2,2,2,2\}}}}="f"
"a":"b1"
"a":"b2"
"a":"c1"
"a":"c2"
"b1":"c3"
"b2":"c3"
"c1":"d1"
"c1":"d2"
"c2":"d2"
"c3":"d3"
"d1":"e"
"d1":"f"
"d2":"e"
"c3":"e"
"b2":@/_1cm/"d1"
}
$$
\end{example}


We skip the proofs of some of the following straightforward results.
Lemma \ref{Lgraph1} is used in \cite{AGP}, while Lemmas \ref{Lgraph3} and \ref{Lgraph4} are used in Section \ref{Sirred}.

\begin{lemma} \label{Lgraph1}
The set $\Omega_{d,0}$ is nonempty if and only if $p \nmid (d+1)$.
If $p \nmid (d+1)$, then $\Omega_{d,0}$ contains one
partition $\{d+2\}$ which is an initial vertex of $G_d$.
If $p \mid (d+1)$, then $\Omega_{d,1}$ consists of $\lceil
(d+1)(p-2)/2p \rceil$ partitions, and every vertex of $G_d$ is
larger than one of these.
\end{lemma}



\begin{lemma} \label{Lgraph3}
If $p=2$, there is a unique maximal partition $\{2, \ldots, 2\}$ in $\Omega_d$ with
length $d/2 +1$.
\end{lemma}

\begin{lemma} \label{Lgraph4}
Let $p \geq 3$.  A partition is maximal if and only if its
entries all equal two or three. Every integer $r+1$ with $(d-1)/3 \leq
r \leq d/2$ occurs exactly once as the length of a maximal
partition. There are $\lfloor d/2 \rfloor - \lceil (d-4)/3 \rceil$
maximal partitions.
There is a unique maximal partition if and only if $d \in
\{1,2,3,5\}$.
\end{lemma}

\begin{proof}
The first statement is true since if $e \geq 4$ then there are $e_1, e_2 \in \ZZ_{>0}$ so that
$e_j \not \equiv 1 \bmod p$ and $e_1+e_2=e$.
For the other statements, let $\vec{E}$ be a maximal partition of $d+2$.
Let $b$ denote the number of the entries of $\vec{E}$ which equal $3$.
Note that $0 \leq b \leq (d+2)/3$.
Let $r+1$ be the length of $\vec{E}$.
Then $d+2=2(r+1)+b$ and $(d+2)/3 \leq r+1 \leq (d+2)/2$.
Any choice of $r+1$ in this range yields a unique choice of $b$ which determines a unique partition $\vec{E}$.
\end{proof}


\begin{remark}
When $p=2$, every path in $G_d$ from the partition $\{d+2\}$ to
the partition $\{2, \ldots, 2\}$ has the same length, which is $d/2$.
When $p=3$, every path in $G_d$ from a minimal to a maximal
vertex has the same length, which is $\lfloor d/3 \rfloor$. This
property does not hold in general for $p \geq 5$.
\end{remark}


\subsection{Artin-Schreier curves} \label{Sascurve}

Here is a review of some basic Artin-Schreier theory.
Let $Y$ be an Artin-Schreier $k$-curve.
Then there is a $\ZZ/p$-cover $\phi:Y \to \PP_k$
with an affine equation of the form $y^p-y=f(x)$ for some non-constant rational function $f(x) \in k(x)$.
At each ramification point, there is a filtration of the inertia group $\ZZ/p$, called the
filtration of higher ramification groups in the lower numbering \cite[IV]{Se:lf}.

Let $\{P_1,\ldots,P_{r+1}\}$ be the set of poles of $f(x)$ on the projective line $\PP_k$.
Let $d_j$ be the order of the pole of $f(x)$ at $P_j$.
One may assume that $p \nmid d_j$ by Artin-Schreier theory.
Then $d_j$ is the {\it lower jump} at $P_j$, i.e., the last index for
which the higher ramification group above $P_j$ is nontrivial.
Let $e_j=d_j+1$.
Then $e_j \geq 2$ and $e_j\not\equiv 1\bmod p$.
The ramification divisor of $\phi$ is $D:=\sum_{j=1}^{r+1}e_j P_j$.

\begin{lemma} \label{Lgenus}
The genus of $Y$ is $g_Y=((\sum_{j=1}^{r+1} e_j) -2)(p-1)/2$.
The $p$-rank of $Y$ is $s_Y=r(p-1)$.
\end{lemma}

\begin{proof}
The first statement follows from the Riemann-Hurwitz formula using \cite[IV, Prop.\ 4]{Se:lf} and
the second from the Deuring-Shafarevich formula \cite[Cor.\ 1.8]{Crew}.
See \cite[Remark 1.4]{Zhu:expsums} or \cite[Section 2]{Z:noextra}
for details.
\end{proof}

\subsection{The $p$-rank of Artin-Schreier curves and partitions}

The Artin-Schreier curves of genus $g=d(p-1)/2$ with $p$-rank $r(p-1)$ are
intimately related to the partition sets $\Omega_{d,r}$ as defined in Section \ref{S3}.

\begin{lemma} \label{Lnonempty}
There exists an Artin-Schreier $k$-curve of genus $g$ with
$p$-rank $r(p-1)$ if and only if $d:=2g/(p-1)$ is a nonnegative integer
and $\Omega_{d,r}$ is nonempty.
\end{lemma}

\begin{proof}
By Lemma \ref{Lgenus}, the existence of an
Artin-Schreier $k$-curve with genus $g=d(p-1)/2$ and $p$-rank $r(p-1)$
is equivalent to the existence of $f(x) \in k(x)$ whose poles have
orders $\{e_1-1, \ldots, e_{r+1}-1\}$ where each $e_j\not\equiv 1\bmod p$
and $\sum_{j=1}^{r+1} e_j=d+2$.
This is equivalent to $\Omega_{d,r}$ being nonempty.
\end{proof}

\begin{example} Let $p=2$.  Let $g \geq 0$ and $0 \leq s \leq g$.
Then $\Omega_{2g,s}$ is non-empty since $2g+2$ can be partitioned into $s+1$ even integers.
Therefore, there exists an Artin-Schreier $k$-curve of genus $g$ and $p$-rank $s$ in
characteristic $2$.
\end{example}

\begin{example}
Let $p \geq 3$.
There exists an Artin-Schreier $k$-curve of genus $g=d(p-1)/2$ with
$p$-rank $0$ if and only if $p \nmid (d+1)$ by Lemma \ref{Lgraph1}.
There exists an ordinary Artin-Schreier $k$-curve (i.e., with $p$-rank $g$)
if and only if $2 \mid d$.
If $2 \nmid d$, the largest $p$-rank which occurs for an Artin-Schreier $k$-curve
of genus $g$ is $s=g-(p-1)/2$ by Lemma \ref{Lgraph4}.
\end{example}

\section{Moduli spaces of Artin-Schreier curves} \label{S4}

Consider fixed parameters $p$, $g=d(p-1)/2$ with $d \geq 1$, and $s=r(p-1)$
with $0 \leq s \leq g$.
In this section, we study the $p$-rank $s$ strata $\AS_{g,s}$ of the moduli space $\AS_g$ of
Artin-Schreier curves of genus $g$.
We show the irreducible components of $\AS_{g,s}$ are in bijection with
the elements of $\Omega_{d,r}$ and find the dimensions of these components.

\subsection{Artin-Schreier covers}

Let $S$ be a $k$-scheme.
An {\it $S$-curve} is a proper flat morphism $Y \to S$ whose
geometric fibres are smooth connected curves.
An {\it Artin-Schreier curve} $Y$ over $S$ is an $S$-curve
for which there exists an (unspecified) inclusion $\iota:\ZZ/p \hookrightarrow \Aut_S(Y)$
such that the quotient $Y/\iota(\ZZ/p)$ is a ruled scheme.
This means that there is an (unspecified) isomorphism between each geometric fibre of $Y/\iota(\ZZ/p)$ and $\PP$.
An {\it Artin-Schreier cover} over $S$
is a $\ZZ/p$-cover $\phi:Y \to \PP_S$.
In other words, it is an Artin-Schreier curve $Y$ over $S$
along with the data of a specified inclusion $\iota:\ZZ/p
\hookrightarrow \Aut_S(Y)$ and a specified isomorphism $Y/\iota(\ZZ/p) \simeq \PP_S$.

Consider the following contravariant functors from the category of
$k$-schemes to sets: $\AS_{g}$ (resp.\ $\ArSc_g$) which associates
to $S$ the set of isomorphism classes of Artin-Schreier curves (resp.\ covers)
over $S$ with genus $g$. As in \cite[Prop.\ 2.7]{Maug},
one can show that there is an algebraic stack
representing $\AS_g$ which we denote again by the symbol $\AS_g$.
Similarly, e.g., \cite[pg.\ 1]{Maug:hur}, there is an algebraic stack representing
$\ArSc_g$ which we denote again by the symbol $\ArSc_g$.
The next lemma is about a natural map from $\ArSc_g$ to $\AS_{g}$.

\begin{lemma} \label{Pforgetful}
Let $g\geq 2$.
There is a morphism $F:\ArSc_g \to \AS_{g}$ and the fibre of $F$
over every geometric point of $\AS_g$ has dimension $3$.
\end{lemma}

\begin{proof}
There is a functorial transformation $\ArSc_g(S) \to \AS_g(S)$
that takes the isomorphism class of a given Artin-Schreier cover $\phi:Y \to \PP_S$ over $S$
to the isomorphism class of the Artin-Schreier curve $Y$ over $S$.
In other words, the transformation is defined by forgetting the inclusion $\iota$
and the isomorphism $Y/\iota(\ZZ/p) \simeq \PP$
(and taking the quotient of the set of inclusions $\iota$ by the action of $\Aut(\ZZ/p)$).
This transformation yields a morphism $F:\ArSc_g \to \AS_{g}$ by Yoneda's lemma.

To prove the second claim, it suffices to work locally in the \'etale topology.
Given an Artin-Schreier $S$-curve $Y$, by \cite[pg.\ 232]{BM}, after an \'etale extension of $S$,
there exists an inclusion $\iota:\ZZ/p \hookrightarrow \Aut_S(Y)$ and an isomorphism $I:Y/\iota(\ZZ/p) \to \PP_S$.
Thus $Y$ is in the image of $F$.
There are only finitely many choices for $\iota$
since $\Aut_S(Y)$ is finite for $g\geq 2$ (e.g., \cite[Theorem 1.11]{DM}).
There is a three-dimensional choice for the isomorphism $I$ since ${\rm dim}({\rm Aut}(\PP_S))=3$.
Thus the fibre of $F$ over $Y$ has dimension three.
\end{proof}

\subsection{The ramification divisor}

This section is about the ramification divisor of a given Artin-Schreier cover.

Let $S$ be a $k$-scheme and let $n \in \ZZ_{> 0}$.
Consider the contravariant functor $\Cdiv_n$, from the category of $k$-schemes to sets,
which associates to $S$ the set of isomorphism classes of 
relative effective Cartier divisors of $\PP_S$ of constant degree $n$.
This functor is represented by a (Hilbert) scheme which we denote also by $\Cdiv_n$.

There is a discrete invariant $\vec{E}$ which induces a natural stratification $\Cdiv_{n,\vec{E}}$ of $\Cdiv_n$.
To see this, suppose $S=\Spec(K)$ where $K$ is a field with ${\rm char}(K)=p$. 
Given $D \in \Cdiv_n(S)$, 
one can associate to $D$ a locally principal effective Weil divisor of $\PP_S$ with degree $n$
by \cite[II, Prop.\ 6.11, Remark 6.11.2]{Hart}.
After a finite flat extension $S' \to S$, 
one can write $D'=D \times_S S'=\sum_{j=1}^{r+1} e_j P_j$ where $e_j \geq 1$ and $\sum_{j=1}^{r+1} e_j=n$
and where $\{P_1, \ldots, P_{r+1}\}$ is a set of distinct horizontal sections of $\PP_{S'}$.
Let $\vec{E}(D)=\{e_1, \ldots, e_{r+1}\}$ with $e_1 \leq \cdots \leq e_{r+1}$.
The partition $\vec{E}(D)$ of $n$ induces a natural stratification $\Cdiv_{n,\vec{E}}$ of $\Cdiv_n$
(where the sections $\{P_1, \ldots, P_{r+1}\}$ associated to $D$ can vary).
Let $\Cdiv^1_n=\cup_{\vec{E} \in \Omega_{n-2}} \Cdiv_{n, \vec{E}}$
(i.e., all $\vec{E}$ for which each $e_j \not \equiv 1 \bmod p$).
For fixed $\vec{E}$, let $H_{\vec{E}} \subset S_{r+1}$ be the subgroup of the symmetric group
generated by all transpositions $(j_1, j_2)$ for which $e_{j_1}=e_{j_2}$.

\begin{lemma} \label{Lbase}
If $\vec{E} \in \Omega_{d,r}$, then
$\Cdiv_{d+2, \vec{E}}$ is irreducible of dimension $r+1$.
\end{lemma}

\begin{proof}
Let $\Delta$ denote the weak diagonal of $(\PP)^{r+1}$,
consisting of $(r+1)$-tuples with at least two coordinates equal.
The quotient of $(\PP)^{r+1}-\Delta$ by the action of $H_{\vec{E}}$
is irreducible with dimension $r+1$.
By the remarks preceding this lemma, 
the spaces $\Cdiv_{d+2, \vec{E}}$ and $[(\PP)^{r+1}-\Delta]/H_{\vec{E}}$
are locally isomorphic for the finite flat topology,
where the isomorphism identifies $D$ with the equivalence class of $(P_1, \ldots, P_{r+1})$.
Thus $\Cdiv_{d+2, \vec{E}}$ is irreducible with dimension $r+1$.
\end{proof}

\begin{proposition} \label{PmapB}
Let $d=2g/(p-1)$.  There is a morphism $B:\ArSc_g \to \Cdiv_{d+2}$ and the
image of $B$ is $\Cdiv^1_{d+2}$.
\end{proposition}

\begin{proof}
Given an Artin-Schreier cover $\phi$ over $S$, consider the closed subscheme $D$ of the 
fixed points under $\iota(\ZZ/p)$; this is a relative Cartier divisor of $\PP_S$ of constant degree $d+2$
by \cite[Lemma 5.2.3, pg.\ 232]{BM}.
The functorial transformation $\ArSc_g(S) \to \Cdiv_{d+2, \vec{E}}(S)$
yields a morphism $B: \ArSc_g \to \Cdiv_{d+2}$ by Yoneda's lemma.

If $D \in \Cdiv_{d+2}(S)$, then there is a restriction on $\vec{E}(D)$.
As before, one can identify a pullback $D'=D \times_S S'$ with
an effective Weil divisor $\sum_{j=1}^{r+1} e_j P_j$ 
where $e_j \geq 1$ are such that $\sum_{j=1}^{r+1} e_j=d+2$ and where
$\{P_1, \ldots, P_{r+1}\}$ is a set of distinct horizontal sections of $\PP_{S'}$.
If $D=B(\phi)$ is in the image of $B$, then $\{P_1, \ldots, P_{r+1}\}$ constitutes the branch
locus of the pullback $\phi'=\phi \times_S S'$ of $\phi$ and $d_j=e_j-1$ is the lower jump of $\phi'$
above the geometric generic point of $P_j$ by \cite[IV, Prop.\ 4]{Se:lf}.
As seen in Section \ref{Sascurve}, $e_j \not \equiv 1 \bmod p$ for $1 \leq j \leq r+1$.
Thus the image of $B$ is contained in $\Cdiv^1_{d+2}$.

Suppose $\vec{E} \in \Omega_{d}$.
To prove that $\Cdiv^1_{d+2, \vec{E}}$ is contained in the image of $B$,
by descent, it suffices to work locally in the finite flat topology.
Given $D=\sum_{j=1}^{r+1} e_j P_j$, consider the divisor $\tilde{D}=\sum_{j=1}^{r+1} (e_j-1)P_j$ of $\PP_{S}$.
There is a non-constant function $f(x)\in \CO(S)(x)$ with ${\rm div}_\infty(f(x))=\tilde{D}$.
Consider the cover $\phi:Y \to \PP_{S}$ given by the affine equation $y^p-y=f(x)$.
Then $\phi$ is an Artin-Schreier cover with ramification divisor $D$ and the fibres of $Y$ have genus $g$
by Lemma \ref{Lgenus}.  Thus $D \in {\rm Im}(B)$.
\end{proof}

Let $\AS_{g,\vec{E}}$ (resp.\ $\ArSc_{g,\vec{E}}$) denote the locally closed reduced subspace of
$\AS_g$ (resp.\ $\ArSc_g$) whose geometric points correspond to Artin-Schreier covers whose ramification divisor
has partition $\vec{E}$.
The morphisms $F$ and $B$ respect the partition $\vec{E}$.
Let $F_{\vec{E}}: \ArSc_{g,\vec{E}} \to
\AS_{g,\vec{E}}$ and $B_{\vec{E}}: \ArSc_{g,\vec{E}} \to
\Cdiv_{d+2, \vec{E}}$ denote the natural restrictions.

\subsection{Artin-Schreier covers with fixed ramification divisor}\label{SfibreB}

In this section, we fix a partition $\vec{E} \in \Omega_{g,r}$
and a divisor $D \in \Cdiv_{d+2, \vec{E}}$ and study the fibre of $B$ over $D$.
Using \cite[Section 5.1]{BM}, we show that this fibre is irreducible and compute its dimension.
We provide some intuition by describing the equations for an Artin-Schreier cover with ramification divisor $D$.

\begin{notation} \label{Ndim}
Let $\vec{E} \in \Omega_{d,r}$ be a fixed partition $\{e_1, \ldots, e_{r+1}\}$ of $d+2$.
Consider a fixed divisor $D$ corresponding to a point of $\Cdiv_{d+2, \vec{E}}$.
Let $\ArSc_{g,D}$ be the fibre of $B_{\vec{E}}: \ArSc_{g, \vec{E}} \to \Cdiv_{d+2, \vec{E}}$ over $D$.
\end{notation}

\begin{notation} \label{Nmoduli}
For $1 \leq j \leq r+1$,
let $t_j=d_j - \lfloor d_j/p \rfloor$ where $d_j=e_j-1$.
Let $N_{\vec{E}}=\sum_{j=1}^{r+1} t_j$.
Let $M_j=(\Aa)^{t_j-1} \times (\Aa-\{0\})$.
Let $M=\times_{j=1}^{r+1} M_j$.
There is an action on $M$ by the subgroup $H_{\vec{E}} \subset S_{r+1}$ generated by all transpositions
$(j_1, j_2)$ for which $d_{j_1}=d_{j_2}$.
Define $M_D= M/H_{\vec{E}}$.
\end{notation}

\begin{proposition} \label{Pdimension}
With notation as in \ref{Ndim} and \ref{Nmoduli},
the fibre $\ArSc_{g,D}$ of $B_{\vec{E}}$ over $D$ is locally isomorphic for the finite flat topology to $M_D$.
Thus $\ArSc_{g,D}$ is irreducible with dimension $N_{\vec{E}}$ over $k$.
\end{proposition}

\begin{proof}
By the definition of $M_D$, the first claim implies the second.
For the first claim,
let $\eta$ denote a labeling of the $r+1$ points in the support of $D$.
Let $\ArSc_{g,D}^{\eta}$ be the contravariant functor which associates to $S$
the set of covers $\phi$ in the fibre $\ArSc_{g,D}(S)$ along with a labeling $\eta$ of the branch locus.
It suffices to show that the moduli space for $\ArSc_{g,D}^{\eta}$ is locally isomorphic to $M$.
This statement can be found in \cite[pg.\ 229, pg.\ 233]{BM}.
\end{proof}

\begin{remark}
In \cite[Cor.\ 2.10]{Ha:mod}, the author constructs an ind-scheme ${\mathcal M}$ which is a
fine moduli space for covers $Y \to \PP$ of $k$-schemes with group $\ZZ/p$ and
branch locus $\{P_1, \ldots, P_{r+1}\}$ (where $Y$ has unbounded genus).
The $k$-points of $\ArSc_{g,D}$ are in bijection with the $k$-points of ${\mathcal M}$ such that $Y$ has genus $g$.
Recall from \cite{Ha:mod} that ${\mathcal M}$ is a direct limit of affine schemes.
This direct limit arises because if $S=\Spec(K)$ where $K$ is not perfect,
then there are non-trivial Artin-Schreier covers over $S$ which become trivial
after a finite flat extension of $S$.  In \cite{Pr:fam}, the author addressed this issue
using a {\it configuration space} whose $k$-points are in bijection with covers defined over $k$.
In Proposition \ref{Pdimension}, we instead followed the approach of \cite[Section 5.1]{BM}.
\end{remark}

\begin{remark}
Proposition \ref{Pdimension} implies that  
a flat base change of $\ArSc_{g, \vec{E}}$ is a ${\mathbb G}_a^n \times {\mathbb G}_m$-bundle over 
$\Cdiv_{d+2, \vec{E}}$ for some $n$.
\end{remark}

\begin{remark}
For the convenience of the reader, we provide some intuition about Proposition \ref{Pdimension}.
Let $S$ be an irreducible affine $k$-scheme.
Suppose $\phi \in \ArSc_{g,D}(S)$ is an Artin-Schreier cover over $S$ with ramification divisor $D$.
Then $\phi$ has an affine equation $y^p-y=f(x)$ for some $f(x) \in \CO(S)(x)$.
The automorphism $\sigma=\iota(1)$ acts via $\sigma(y)=y+z$ for some $z \in (\ZZ/p)^*$.
Two such covers $\phi_1:y^p-y=f_1(x)$ and $\phi_2:y^p-y=f_2(x)$ are isomorphic if and only if
$f_2(x)=(z_2/z_1) f_1(x) +\delta^p-\delta$ for some $\delta \in \CO(S)(x)$, see e.g., \cite[Lemma 2.1.5]{Pr:fam}.
After possibly changing $f(x)$, one can suppose $z=1$.

The cover $\phi$ is in {\it standard form} if
$p \nmid i$ for any monomial $c_ix^i$ in $f(x)$ whose
coefficient $c_i$ is generically non-nilpotent.
Given an Artin-Schreier cover $\phi$, after a finite flat extension $S' \to S$,
then $\phi \times_S S'$ has an affine equation in standard form.
To prove this, one uses an \'etale cover $S'' \to S$ with equation $a^p-a=c_0$ to remove a constant
coefficient $c_0 \in \CO(S)$ from $f(x)$.  If $f(x)$ contains a monomial $cx^{pw}$ with $w \in \ZZ_{>0}$,
one uses a purely inseparable cover $S' \to S''$ with equation $b^p=c$ to replace $cx^{pw}$ with the monomial $bx^{w}$.
These transformations are uniquely determined and do not change the isomorphism class of $\phi \times_S S'$.

Suppose $D=\sum_{j+1}^{r+1} e_j P_j$ where $\{P_1, \ldots P_{r+1}\}$ is a fixed set of distinct horizontal sections of $\PP_S$.
If $\phi$ has ramification divisor $D$, then $f(x)$ has a partial fraction decomposition
$f(x)=\sum_{j=1}^{r+1} g_j(x)$ where $g_j(x)\in (x-P_j)^{-1}\CO(S)[(x-P_j)^{-1}]$
is a polynomial of degree $d_j$ in the variable $(x-P_j)^{-1}$ with no constant term.
(If $P_j=P_\infty$, let $(x-P_j)^{-1}$ denote $x$ for consistency of notation.)
If $\phi$ is in standard form, one can write $g_j(x)=\sum_{i=1}^{d_j} c_{i,j}(x-P_j)^{-i}$
where $c_{i,j}=0$ if $p \mid i$ and $c_{d_j,j}$ is never zero.
The isomorphism between $M$ and $\ArSc_{g,D}^{\eta}$ in the finite flat topology
identifies $(\times_{j=1}^{r+1} \times_{i=1, \ p \nmid i}^{d_j} c_{i,j})$
with the isomorphism class of the Artin-Schreier cover
$y^p-y=\sum_{j=1}^{r+1} g_j(x)$ (with the implicit labeling of $\{P_1, \ldots, P_{r+1}\}$).
\end{remark}

\subsection{Irreducible components of the $p$-rank strata} \label{Scover}

Recall that $g=d(p-1)/2$ with $d \geq 1$ and $d$ even if $p=2$ and $s=r(p-1)$ with $0 \leq s \leq g$.
The $p$-rank induces a stratification of $\AS_g$ (resp.\ $\ArSc_g$).
Let $\AS_{g,s}$ (resp.\ $\ArSc_{g,s}$) denote the locally closed reduced subspace of
$\AS_g$ (resp.\ $\ArSc_g$) whose geometric points have $p$-rank $s$.

\begin{theorem} \label{Tascov}
The irreducible components of $\ArSc_{g,s}$ are the strata $\ArSc_{g,\vec{E}}$ with $\vec{E} \in \Omega_{d,r}$.
If $\vec{E}=\{e_1, \ldots, e_{r+1}\}$, then the dimension over $k$ of the irreducible component $\ArSc_{g,\vec{E}}$
is $d+2-\sum_{j=1}^{r+1} \lfloor (e_j-1)/p \rfloor$.
\end{theorem}

\begin{proof}
The image of $\ArSc_{g,s}$ under $B$ is the union of the strata $\Cdiv_{d+2, \vec{E}}$ of
$\Cdiv^1_{d+2}$ with $r(\vec{E})=r$ by Proposition \ref{PmapB}.
The stratum $\Cdiv_{d+2, \vec{E}}$ is irreducible of dimension $r+1$ by Lemma \ref{Lbase}.

For $\vec{E} \in \Omega_{d,r}$, consider the
morphism $B_{\vec{E}}:\ArSc_{g,\vec{E}} \to \Cdiv_{d+2, \vec{E}}$.
The fibre of $B_{\vec{E}}$ over a fixed
divisor $D$ is irreducible by Proposition \ref{Pdimension}.
By Zariski's main theorem, $\ArSc_{g, \vec{E}}$ is irreducible since $B_{\vec{E}}$
has irreducible fibres and image.
Thus the irreducible components of $\ArSc_{g,s}$ are the strata $\ArSc_{g,\vec{E}}$ with $\vec{E} \in \Omega_{d,r}$.

The morphism $B_{\vec{E}}$ is flat since all its fibres are isomorphic.
Thus the dimension of $\ArSc_{g,\vec{E}}$ is the sum of the dimensions
of $\Cdiv_{d+2, \vec{E}}$ and of the fibres of $B_{\vec{E}}$.
This equals $r+1+\sum_{j=1}^{r+1}(d_j - \lfloor d_j/p \rfloor)$ by Lemma \ref{Lbase} and Proposition \ref{Pdimension}.
This simplifies to $d+2-\sum_{j=1}^{r+1} \lfloor (e_j-1)/p \rfloor$.
\end{proof}

Theorem \ref{T1} in the introduction follows immediately from the next corollary.

\begin{corollary} \label{Maincor}
The irreducible components of $\AS_{g,s}$ are the strata $\AS_{g, \vec{E}}$ with $\vec{E} \in \Omega_{d,r}$.
If $\vec{E}=\{e_1, \ldots, e_{r+1}\}$, then the dimension $d_{\vec{E}}$ over $k$ of the irreducible component $\AS_{g,\vec{E}}$
is $d-1-\sum_{j=1}^{r+1} \lfloor (e_j-1)/p \rfloor$.
\end{corollary}

\begin{proof}
Let $W$ be an irreducible component of $\AS_{g,s}$.
By Lemma \ref{Pforgetful}, $F^{-1}(W)$ is a union of irreducible components of $\ArSc_{g,s}$.
By Theorem \ref{Tascov}, these are indexed by partitions $\vec{E} \in \Omega_{d,r}$.
The morphism $F$ respects the partition $\vec{E}$.
In other words, given an Artin-Schreier curve $Y$,
every Artin-Schreier cover $\phi:Y \to \PP_k$ has the same partition.
Thus there is a unique partition occurring for points in $F^{-1}(W)$, and so $F^{-1}(W)$ is irreducible.
So the irreducible components of $\AS_{g,s}$ are the strata $\AS_{g, \vec{E}}$ with $\vec{E} \in \Omega_{d,r}$.
The second statement follows by Lemma \ref{Pforgetful} for $g \geq 2$ since $\dime(W)=\dime(F^{-1}(W))-3$
and by direct computation for $g=1$.
\end{proof}

\begin{example}
Let $p=3$ and $g=10$.  Here are the dimensions $d_{\vec{E}}$ of the irreducible components
of $\AS_{10,s}$.
\[\begin{array}{|l|l|}
\hline
s & {\rm dimension}\\
\hline
0 & d_{\{12\}}=6 \\
\hline
2 & d_{\{3,9\}}=7, \ d_{\{6,6\}}=7 \\
\hline
4 & d_{\{2,2,8\}}=7, \ d_{\{2,5,5\}}= 7, \  d_{\{3,3,6\}}=8  \\
\hline
6 & d_{\{2,2,2,6\}}=8, \ d_{\{2,2,3,5\}}=8, \ d_{\{3,3,3,3\}}=9 \\
\hline
8 & d_{\{2,2,2,3,3\}}=9 \\
\hline
10 & d_{\{2,2,2,2,2,2\}}=9 \\
\hline
\end{array}
\]
\end{example}

The next corollary shows that the image of $\AS_g$ under the Torelli morphism is not in general position
relative to the $p$-rank stratification of $\CA_g$ when $p \geq 3$.

\begin{corollary} \label{Cnottrans}
If $p \geq 3$, then ${\rm codim}(\AS_{g,s},\AS_g) < g-s$.
\end{corollary}

\begin{proof}
Let $d=2g/(p-1)$ and $r=s/(p-1)$.
Let $\epsilon = \min\sum_{j=1}^{r+1}\lfloor (e_j-1)/p\rfloor$
where the minimum ranges over all partitions $\{e_1,\ldots,e_{r+1}\}$ with fixed sum $d+2$.
By Corollary \ref{Maincor}, ${\rm codim}(\AS_{g,s},\AS_g)=\epsilon$.
Since $\lfloor (e_j-1)/p\rfloor \leq (e_j-2)/p$, one sees that $\epsilon \leq (d-2r)/p=2(g-s)/p(p-1)$.
Thus $\epsilon < g-s$ if $p \geq 3$.
\end{proof}

\subsection{Irreducibility of the Artin-Schreier locus} \label{Sirred}

As an application of Theorem \ref{T1}, we determine all pairs $(p,g)$ for which $\AS_g$ is irreducible.

\paragraph{Corollary \ref{C2}
\it
The moduli space $\AS_g$ is irreducible in exactly the following cases:
(i) $p=2$; or (ii) $g=0$ or $g = (p-1)/2$; or (iii) $p=3$ and $g=2,3,5$.}

\begin{proof}
Let $d=2g/(p-1)$.
Recall that $d_{\vec{E}}$ is the dimension of $\AS_{g, \vec{E}}$.
The first claim is that there is a bijection between irreducible components of $\AS_g$
and partitions $\vec{E} \in \Omega_{d,r}$ so that $d_{\vec{E}}=d-1$.
To see this, note that \cite[Cor.\ 3.16]{Maug} implies that every irreducible component of $\AS_g$ has dimension $d-1$.
If $\Gamma$ is an irreducible component of $\AS_g$, then there is a partition $\vec{E} \in \Omega_d$
and an open subset $U \subset \Gamma$ so that $U \subset \AS_{g, \vec{E}}$.
Then $d_{\vec{E}}=\dime(\Gamma)=d-1$.
Conversely, suppose $d_{\vec{E}}=d-1$ for some $\vec{E} \in \Omega_{d}$.
Then the irreducible space $\AS_{g,\vec{E}}$ is open in a unique irreducible component $\Gamma$ of $\AS_g$.

Thus, $\AS_g$ is irreducible if and only if there is exactly one partition
$\vec{E} \in \Omega_d$ with dimension $d_{\vec{E}}=d-1$.
Write $\vec{E}=\{e_1, \ldots, e_{r+1}\}$.
By Theorem \ref{T1}, $d_{\vec{E}}=d-1$ if and only if $e_j < p+1$ for $1 \leq j \leq r+1$.

If $p=2$, only one partition satisfies the condition $e_j < 3$ for each $j$,
namely the partition $\{2, \ldots, 2\}$, Lemma \ref{Lgraph3}.  Thus $\CA_g$ is irreducible for all $g$ when $p=2$.

For arbitrary $p$, if $g=0$ (resp.\ $g=(p-1)/2$) then $d=0$ (resp.\ $d=1$), and there is only one partition
satisfying $e_j < p+1$, namely the partition $\{2\}$ (resp.\ $\{3\}$).
Thus $\AS_g$ is irreducible in these cases.

If $p=3$ and $d=2$ (resp.\ $3$, $5$), only one partition satisfies the
condition $e_j < 4$, namely $\{2,2\}$, (resp.\ $\{2,3\}$, $\{2,2,3\}$).
Thus $\AS_g$ is irreducible in these cases.

Suppose $p \geq 3$ and $d \geq 2$ and that $\AS_g$ is irreducible.
If $\vec{E}$ is a maximal partition, then its entries satisfy $e_j \leq 3 < p+1$.
Thus $\Omega_d$ has a unique maximal partition.
By Lemma \ref{Lgraph4}, this implies $d \in \{2,3,5\}$.
If $p \geq 5$, then there are at least two partitions satisfying $e_j < p+1$:
for example, $\{4\}$ and $\{2,2\}$ when $d=2$; $\{5\}$ and $\{2,3\}$ when $d=3$;
$\{2,5\}$ and $\{2,2,3\}$ when $d=5$.
This is a contradiction and so
cases (i)-(iii) are the only cases when $\AS_g$ is irreducible.
\end{proof}

\subsection{Hyperelliptic curves in characteristic $2$} \label{Shyper}

Let $\CH_g$ be the moduli space of hyperelliptic $k$-curves of genus $g$.
Let $\CH_{g,s}$ denote the locally closed reduced subspace of $\CH_g$
parametrizing hyperelliptic $k$-curves of genus $g$ with $p$-rank $s$.
When $p=2$, $\CH_g$ is the same as $\AS_g$.
This yields the following result.

\paragraph{Corollary \ref{C1}}
{\it Let $p=2$.  The irreducible components of $\CH_{g,s}$ are in bijection
with partitions of $g+1$ into $s+1$ positive integers.
Every component has dimension $g-1+s$ over $k$.}

\begin{proof}
By Corollary \ref{Maincor}, the irreducible components of $\CH_{g,s}$ are in bijection with
the partitions of $d+2=2g+2$ into $s+1$ even positive integers,
which are in bijection with the partitions of $g+1$ into $s+1$
positive integers. The dimension of the irreducible component for
$\vec{E}=\{e_1, \ldots, e_{s+1}\}$ is
$(d-1)-\sum_{j=1}^{s+1} \lfloor (e_j-1)/2 \rfloor$.
This simplifies to $g-1+s$ since
$e_j$ is even and $\lfloor (e_j-1)/2 \rfloor=e_j/2 -1$.
\end{proof}

\section{Deformation results and open questions} \label{Sdeform}

In this section, we give some results on how the irreducible components of
$\AS_{g,s}$ (with varying $s$) fit together within $\AS_g$.
This involves deformations of wildly ramified covers with non-constant branch locus.

\subsection{A deformation result for wildly ramified covers} \label{Ssubdeform}

The main result of this section is that, under the obvious necessary conditions,
the $p$-rank of an Artin-Schreier curve can be increased by exactly $p-1$ in a flat deformation.
Let $S=\Spec(k[[t]])$ and let $s$ be the closed point of $S$.

\begin{proposition} \label{Pdeformp}
Suppose $p \mid e_1$ or $p \mid e_2$.
Suppose $\psi_\circ$ is an Artin-Schreier cover over $k$,
branched at a point $b$ with lower jump $e_1+e_2-1$.
Then there exists an Artin-Schreier cover $\psi_S$ over $S$
whose special fibre is isomorphic to $\psi_\circ$, whose generic fibre is branched at
two points that specialize to $b$ and which have lower jumps $e_1-1$ and $e_2-1$, and whose
ramification divisor is otherwise constant.
\end{proposition}

\begin{proof}
Let $e=e_1+e_2$.  By hypothesis, $p \nmid (e-1)$.  Without loss of generality, suppose $p \mid e_1$.

Consider the Artin-Schreier cover $\psi_\circ:Y_\circ \to Z_\circ$ which is wildly ramified at the point
$y_\circ \in Y_\circ$ above $b$ where it has lower jump $e-1$.
Let $\hat{\psi}_\circ:\hat{Y}_{\circ} \to \hat{Z}_{\circ}$ be the germ of
$\psi_\circ$ at $y_\circ$.  It is an Artin-Schreier cover of germs of curves.
Using formal patching, see e.g., \cite[Prop.\ 2.7]{Ha:mod} or \cite[Thm.\ 3.3.4]{BM},
deformations of $\psi_\circ$ can be constructed locally via deformations of $\hat{\psi}_\circ$.
With this technique, one can suppose that the
deformation of $\psi_\circ$, and thus the ramification divisor, is constant away from $b$.

Now $\hat{Z}_\circ \simeq \Spec(k[[x^{-1}]])$.
After a change of variables, one can suppose that
the restriction of $\hat{\psi}_\circ$ to $\Spec(k((x^{-1})))$ has equation $y^p-y=x^{e-1}$.

Consider the deformation $\hat{\psi}_S$ of $\hat{\psi}_\circ$ over $S=\Spec(k[[t]])$ given by
the normal extension of $\Spec(k[[x^{-1},t]])$ with the following affine equation:
$$y^p-y=x^{e-1}/(1-xt)^{e_1}.$$
On the special fibre, when $t=0$, then $\hat{\psi}_s$ is isomorphic to $\hat{\psi}_\circ$.
On the generic fibre, when $t \not = 0$, then $\hat{\psi}_{S-s}$ is branched above
$x^{-1}=0$ and above $x^{-1}=t$.
Let $F(x)=x^{e-1}/(1-xt)^{e_1}$.
The order of the pole of $F(x)$ at $x^{-1}=0$ is $e-1-e_1=e_2-1$, which is prime-to-$p$ by hypothesis.
Thus the lower jump above $x^{-1}=0$ is $e_2-1$.

To compute the lower jump above $x^{-1}=t$, one can expand $F(x)$ around $1/t$:
$$F(x)=(-1)^{e_1}t^{-(e+e_1-1)}(x-1/t)^{-e_1} + (e-1)(-1)^{e_1}t^{-(e+e_1-2)}(x-1/t)^{-e_1+1} + \ldots.$$
After a finite inseparable extension of $k((t))$ with equation $t_1^p=t$, the leading term of $F(x)$ is a $p$th power.
The second term of $F(x)$ is non-zero since $p \nmid (e-1)$ and thus it
becomes the leading term of the affine equation in standard form for $\hat{\psi}_{S-s}$.
Thus the lower jump above $x^{-1}=t$ is $e_1-1$.
Thus the cover $\hat{\psi}_{S-s}$ is branched at two points that specialize to $b$
and which have lower jumps $e_1-1$ and $e_2-1$.
By Lemma \ref{Lgenus} and \cite[Lemma IV.2.3]{OSS}, the deformation $\hat{\psi}_S$ of $\hat{\psi}_\circ$ over $S$ is smooth.
\end{proof}

The next result shows that, under a mild necessary condition, the $p$-rank of an Artin-Schreier curve
can be increased by exactly $p-1$ under a flat deformation.
In particular, an Artin-Schreier curve of genus $g \geq p(p-1)/2$ and $p$-rank $0$ can be deformed
to an Artin-Schreier curve of genus $g$ and $p$-rank $p-1$.

\begin{proposition} \label{Pdeformup}
Suppose that $Y_\circ$ is an Artin-Schreier $k$-curve of genus $g$ and $p$-rank $r(p-1)$.
Suppose there is a ramified point of $Y_\circ$ under the $\ZZ/p$-action whose
lower jump $d$ satisfies $d \geq p+1$.
Then there exists an Artin-Schreier curve $Y_S$ over $S$
whose special fibre is isomorphic to $Y_\circ$ and whose generic fibre has genus $g$ and $p$-rank $(r+1)(p-1)$.
\end{proposition}

\begin{proof}
Let $e_1=p$ and $e_2 = d+1-p$.
By hypothesis, there is an Artin-Schreier cover $\psi_\circ:Y_\circ \to \PP_k$,
branched at $r+1$ points, including one point $b$ with lower jump $e_1+e_2-1$.
The result is then immediate from Proposition \ref{Pdeformp},
because the generic fibre of $\psi_S$ is branched at $r+2$ points.
\end{proof}

\subsection{Preliminary closure results} \label{Sprelimclosure}

In this section, we show that the combinatorial data in the graph $G_d$ gives partial information about
how the irreducible components of $\AS_{g,s}$ (with varying $s$) fit together in $\AS_g$.
In fact, we will see in Section \ref{Sgraphdeform} that
the graph $G_d$ gives complete information about this question when $p=2$.

For $i=1,2$, consider a partition $\vec{E}_i \in \Omega_{d,r_i}$.
Let $s_i=r_i(p-1)$. Let $\Gamma_{\vec{E_i}}:=\AS_{g, \vec{E}_i}$
be the irreducible component of $\AS_{g,s_i}$ corresponding to
$\vec{E}_i$ as defined below Proposition \ref{PmapB}. There is a
partial ordering $\prec$ on $\Omega_d$ from Section \ref{S3}.

\begin{lemma} \label{Pdeformgraph}
If $\Gamma_{\vec{E}_1}$ is in the closure of $\Gamma_{\vec{E}_2}$
in $\AS_g$, then $\vec{E}_1 \prec \vec{E}_2$.
\end{lemma}

\begin{proof}
Let $S=\Spec(k[[t]])$ and consider an Artin-Schreier cover
$\phi_S$ so that the generic fibre yields a $k((t))$-point of
$\Gamma_{\vec{E}_1}$ and the special fibre yields a $k$-point of $\Gamma_{\vec{E}_2}$. This is
only possible if the branch points of $\phi_S$ coalesce when
$t=0$. Since $B(\phi_S)$ is a relative Cartier divisor of constant degree,
the entries of the partition sum together under specialization
and the partition decreases in size.
\end{proof}

The next example and lemma show that the condition $\vec{E}_1
\prec \vec{E}_2$ is frequently not sufficient for
$\Gamma_{\vec{E}_1}$ to be in the closure of $\Gamma_{\vec{E}_2}$
in $\AS_g$ when $p \geq 5$.

\begin{example}
Let $p=5$ and $g=4$ and consider $\vec{E}_1=\{4\}$ and
$\vec{E}_2=\{2,2\}$. Then $\Gamma_{\vec{E}_1}$ and
$\Gamma_{\vec{E}_2}$ are both components of $\AS_{4}$ with
dimension one. Although $\vec{E}_1 \prec \vec{E_2}$, at most a
zero-dimensional subvariety of $\Gamma_{\vec{E}_1}$ can be in the
closure of $\Gamma_{\vec{E}_2}$. In fact, $\Gamma_{\vec{E}_1}$ is
the supersingular family parametrized by $y^5-y=x^3+cx^2$; while
$\Gamma_{\vec{E}_2}$ is the ordinary family parametrized by
$y^5-y=x+c/x$.
\end{example}

For $a \in \ZZ_{>0}$, let $\overline{a}$ be the integer so that
$\overline{a} \equiv a \bmod p$ and $0 \leq \overline{a} < p$.

\begin{lemma} \label{Ldimcomp}
Suppose $\vec{E}_1 \prec \vec{E}_2$ with an edge from $\vec{E}_1$
to $\vec{E}_2$.
\begin{enumerate}
\item If the edge is of the form $\{e\} \mapsto \{e_1,e_2\}$ with $2 <  \overline{e}_1 +\overline{e}_2 \leq p$,
then $\dime_k(\Gamma_{\vec{E}_1})=\dime_k(\Gamma_{\vec{E}_2})$ and $\Gamma_{\vec{E}_1}$
is not in the closure of $\Gamma_{\vec{E}_2}$ in $\AS_g$.
\item In all other cases, $\dime_k(\Gamma_{\vec{E}_1})=\dime_k(\Gamma_{\vec{E}_2})-1$.
\end{enumerate}
\end{lemma}

\begin{proof}
The dimension comparison follows from Theorem \ref{T1}.
If $\dime_k(\Gamma_{\vec{E}_1})=\dime_k(\Gamma_{\vec{E}_2})$, then $\Gamma_{\vec{E}_1}$ is
not in the closure of $\Gamma_{\vec{E}_2}$
since $\AS_g$ is separated.
\end{proof}

\subsection{Closure of the $p$-rank strata} \label{Sgraphdeform}

The main result of this section is Theorem \ref{Tclosureagain} which states that
every irreducible component of $\AS_{g,s}$ satisfying an obvious necessary condition
is contained in the closure of $\AS_{g,s+(p-1)}$ in $\AS_g$.
In the case that $p=2$, Corollary \ref{deform2} strengthens this result.

Recall, for $i=1,2$, that $\Gamma_{\vec{E_i}}:=\AS_{g, \vec{E}_i}$
is the irreducible component of $\AS_{g,s_i}$ corresponding to $\vec{E}_i \in \Omega_{d, r_i}$,
where $s_i=r_i(p-1)$.
There are some earlier results about when $\Gamma_{\vec{E}_1}$ is in the closure of $\Gamma_{\vec{E}_2}$.
For example, \cite[Thm.\ 6.5.1]{M:thesis} implies that $\Gamma_{\vec{E}_1}$ is in the
closure of $\Gamma_{\vec{E}_2}$ for an edge of the form
$\{2p-\ell+1\} \mapsto \{p, p-\ell+1\}$ as long as $\ell \mid (p-1)$.
Here is another such result.

\begin{proposition} \label{Pboundaryp}
Let $\vec{E}_1 \prec \vec{E}_2$ with an edge of the form $\{e\}
\mapsto \{e_1,e_2\}$ from $\vec{E}_1$ to $\vec{E}_2$. If $p \mid
e_1$ or $p \mid e_2$, then $\Gamma_{\vec{E}_1}$ is in the closure
of $\Gamma_{\vec{E}_2}$ in $\AS_g$.
\end{proposition}

In other words, under the hypothesis of Proposition \ref{Pboundaryp},
if $Y_\circ$ is an Artin-Schreier curve with partition $\vec{E}_1$ over $k$,
then there exists an Artin-Schreier curve $Y_S$ over $S=\Spec(k[[t]])$
whose special fibre is isomorphic to $Y_\circ$ and whose generic fibre has partition $\vec{E}_2$.

\begin{proof}
For $i=1,2$, let $\Gamma_{\vec{E}_i}=\AS_{g, \vec{E}}$.
Let $Y_\circ$ be the Artin-Schreier curve corresponding to a $k$-point of $\Gamma_{\vec{E}_1}$.
There exists an Artin-Schreier cover $\phi_\circ:Y_\circ \to \PP_k$ over $k$.
The element $e$ in the partition $\vec{E}_1$ determines a branch point $b \in \PP_k$
so that the lower jump of $\phi_\circ$ above $b$ is $e-1$.

Let $S=\Spec(k[[t]])$.
By Proposition \ref{Pdeformp}, there exists an Artin-Schreier cover $\phi_S$ over $S$
whose special fibre is isomorphic to $\phi_\circ$ and whose generic fibre is branched at
two points that specialize to $b$ and that have lower jumps $e_1-1$ and $e_2-1$.
Furthermore, the ramification divisor is otherwise constant.
Thus the generic fibre of $\phi_S$ has partition $\vec{E}_2$.
Thus $\Gamma_{\vec{E}_1}$ is in the closure of $\Gamma_{\vec{E}_2}$ in $\AS_g$.
\end{proof}

\begin{theorem} \label{Tclosureagain}
Suppose $0 \leq s \leq g-(p-1)$.
If $\eta$ is an irreducible component of $\AS_{g,s}$ which is not
open and dense in an irreducible component of $\AS_g$,
then $\eta$ is in the closure of $\AS_{g,s+(p-1)}$ in $\AS_g$.
\end{theorem}

\begin{proof}
The condition that $\eta$ is not open and dense in an irreducible component of $\AS_g$
implies that $\dime(\eta)<d-1$ \cite[Cor.\ 3.16]{Maug}.
By Theorem \ref{Maincor}, $\eta=\AS_{g, \vec{E}}$ for some partition $\vec{E} \in \Omega_{d,r}$
containing an entry $e \geq p+2$.
The result is then immediate from Proposition \ref{Pboundaryp}, letting $e_1=p$ and $e_2=e-p$.
\end{proof}

The next corollary shows that the graph $G_d$ gives a complete combinatorial description
of how the irreducible components of $\AS_{g,s}$ fit together in $\AS_g$ when $p=2$.
This result is used in \cite{AGP}.

\begin{corollary} \label{deform2}
Suppose $p=2$. Then $\Gamma_{\vec{E}_1}$ is in the closure of
$\Gamma_{\vec{E}_2}$ in $\AS_g$ if and only if $\vec{E}_1 \prec \vec{E}_2$.
Thus, if $0 \leq s < g$, then every component of $\CH_{g,s}$ is in the closure of $\CH_{g,s+1}$ in $\CH_g$.
\end{corollary}

\begin{proof}
Lemma \ref{Pdeformgraph} implies the forward direction.
For the converse, one reduces to the case that there is an edge from $\vec{E}_1$ to $\vec{E}_2$.
Since $p=2$, the edge has the form $\{e\} \mapsto \{e_1,e_2\}$ where $e_1$ and $e_2$ are even.
Then Proposition \ref{Pboundaryp} applies.
\end{proof}

\subsection{Open questions}

An answer to the following question would help determine whether $\AS_g$ is connected.

\paragraph{Question 1:} \label{Q1}
What are necessary and sufficient conditions on the edge $\{e\} \mapsto \{e_1,e_2\}$
or the edge $\{e\} \mapsto \{e_1,e_2,e_3\}$ for $\Gamma_{\vec{E}_1}$ to be in the closure of $\Gamma_{\vec{E}_2}$ in $AS_g$?

\begin{remark}
By Proposition \ref{Pboundaryp}, a sufficient condition for an affirmative answer to Question \ref{Q1}
for the edge $\{e\} \mapsto \{e_1,e_2\}$ is that $p \mid e_1 e_2$.
Here is a heuristic why this condition may also be necessary.
Suppose that $K$ is a field of characteristic $0$. If $\Phi:Y \to \PP_K$ is a $\ZZ/p$-Galois cover and $y \in Y$ is a ramification point,
then the identification of $\ZZ/p$ with ${\rm Gal}(\Phi)$ allows one to define a
{\it canonical generator} $g_y \in \ZZ/p$ of the inertia group at $y$, see e.g., \cite[Section 2.2.1]{V:book}.
The inertia type of $\Phi$ is the multi-set $\{g_y \}$ for all ramification points $y$ of $\Phi$.

Now, if $\phi_1$ is an Artin-Schreier cover (in characteristic $p$) with partition $\{e\}$,
then $\phi_1$ can be lifted to a $\ZZ/p$-cover of the projective line over a field of characteristic $0$,
and the inertia type of this lifting is the multi-set of length $e$ of the form $\{1, \ldots, 1, 1-e\}$ \cite[Ex.\ 3.3.1]{GM:order}.
Similarly, if $\phi_2$ is an Artin-Schreier cover with partition $\{e_1, e_2\}$,
then the inertia type of the lifting is the multi-set of length $e_1+e_2=e$ of the form
$\{1-e_1, 1, \ldots, 1, 1-e_2\}$. So the inertia types of the liftings are the same if and only if
either $1-e_1 \equiv 1 \bmod p$ or $1-e_2 \equiv 1 \bmod p$,
in other words, if and only if $p \mid e_1 e_2$.
\end{remark}

\paragraph{Question 2:}
Let $\vec{E} \in \Omega_{d,r}$.
What Newton polygons occur for points of $\AS_{g,\vec{E}}$?

\smallskip

When $p \gg d$, the Newton polygon occurring for the generic point of
$\AS_{g,\vec{E}}$ is found in \cite{Zhu:expsums}.
Its limit as $p \to \infty$
has slopes $0$ and $1$ occurring with multiplicity $r(p-1)$ and
slopes $\{\frac{1}{e_j-1},\ldots,\frac{e_j-2}{e_j-1}\}$ with
multiplicity $p-1$ for each $1 \leq j \leq r+1$.
See \cite{blache08}, \cite{blache09} for recent results on this question.

\smallskip

\paragraph{Question 3}
If $p\geq 3$ and $g> s \geq 0$, is every component of $\CH_{g,s}$ in the closure of $\CH_{g,s+1}$?

\smallskip

An answer to Question 3 would give more information about the
geometry of the $p$-rank stratification of $\CH_g$, thus generalizing Corollary \ref{deform2}.
In \cite[Cor.\ 3.15]{AP:prh}, the authors prove a related result: if $p \geq 3$ and $0 \leq s' < s \leq g$,
then for each irreducible component $S$ of $\CH_{g,s}$, there exists an irreducible component $T_{s'}$ of
$\CH_{g,s'}$ such that $\bar S$ contains $T_{s'}$.




\bibliographystyle{plain}
\bibliography{paper34clean}

\begin{thebibliography}{10}

\bibitem{AP:prh}
J.~Achter and R.~Pries.
\newblock The $p$-rank strata of the moduli space of hyperelliptic curves.
\newblock arXiv:0902.4637.

\bibitem{AGP}
Jeffrey~D. Achter, Darren Glass, and Rachel Pries.
\newblock Curves of given {$p$}-rank with trivial automorphism group.
\newblock {\em Michigan Math. J.}, 56(3):583--592, 2008.
\newblock arXiv:0708.2199.

\bibitem{BM}
J.~Bertin and A.~M{\'e}zard.
\newblock D\'eformations formelles des rev\^etements sauvagement ramifi\'es de
  courbes alg\'ebriques.
\newblock {\em Invent. Math.}, 141(1):195--238, 2000.

\bibitem{blache09}
R.~Blache.
\newblock First vertices for generic newton polygons, and {$p$}-cyclic
  coverings of the projective line.
\newblock arXiv:0912.2051.

\bibitem{blache08}
R.~Blache.
\newblock {$p$}-density, exponential sums and {A}rtin-{S}chreier curves.
\newblock arXiv:0812.3382.

\bibitem{Crew}
R.~Crew.
\newblock \'{E}tale {$p$}-covers in characteristic {$p$}.
\newblock {\em Compositio Math.}, 52(1):31--45, 1984.

\bibitem{DM}
P.~Deligne and D.~Mumford.
\newblock The irreducibility of the space of curves of given genus.
\newblock {\em Inst. Hautes \'Etudes Sci. Publ. Math. No.}, 36:75--109, 1969.

\bibitem{FVdG:complete}
C.~Faber and G.~van~der Geer.
\newblock Complete subvarieties of moduli spaces and the {P}rym map.
\newblock {\em J. Reine Angew. Math.}, 573:117--137, 2004.

\bibitem{GP:05}
D.~Glass and R.~Pries.
\newblock Hyperelliptic curves with prescribed {$p$}-torsion.
\newblock {\em Manuscripta Math.}, 117(3):299--317, 2005.

\bibitem{GM:order}
Barry Green and Michel Matignon.
\newblock Order {$p$} automorphisms of the open disc of a {$p$}-adic field.
\newblock {\em J. Amer. Math. Soc.}, 12(1):269--303, 1999.

\bibitem{Ha:mod}
D.~Harbater.
\newblock Moduli of $p$-covers of curves.
\newblock {\em Comm. Algebra}, 8(12):1095--1122, 1980.

\bibitem{Hart}
Robin Hartshorne.
\newblock {\em Algebraic geometry}.
\newblock Springer-Verlag, New York, 1977.
\newblock Graduate Texts in Mathematics, No. 52.

\bibitem{Lonsted}
K.~L{\o}nsted.
\newblock The hyperelliptic locus with special reference to characteristic two.
\newblock {\em Math. Ann.}, 222(1):55--61, 1976.

\bibitem{Maug:hur}
S.~Maugeais.
\newblock On a compactification of a {H}urwitz space in the wild case.
\newblock math.AG/0509118.

\bibitem{Maug}
S.~Maugeais.
\newblock Quelques r\'esultats sur les d\'eformations \'equivariantes des
  courbes stables.
\newblock {\em Manuscripta Math.}, 120(1):53--82, 2006.

\bibitem{M:thesis}
A.~M\'ezard.
\newblock Quelques probl\`emes de d\'eformations en caract\'eristique mixte.
\newblock th\`ese de doctorat de math\'ematiques de l'universit\'e Joseph
  Fourier.

\bibitem{O:purity}
F.~Oort.
\newblock Subvarieties of moduli spaces.
\newblock {\em Invent. Math.}, 24:95--119, 1974.

\bibitem{Pr:fam}
R.~Pries.
\newblock Families of wildly ramified covers of curves.
\newblock {\em Amer. J. Math.}, 124(4):737--768, 2002.

\bibitem{SZ:02}
J.~Scholten and H.~J. Zhu.
\newblock Hyperelliptic curves in characteristic 2.
\newblock {\em Int. Math. Res. Not.}, (17):905--917, 2002.

\bibitem{OSS}
T.~Sekiguchi, F.~Oort, and N.~Suwa.
\newblock On the deformation of {A}rtin-{S}chreier to {K}ummer.
\newblock {\em Ann. Sci. \'Ecole Norm. Sup. (4)}, 22(3):345--375, 1989.

\bibitem{Se:lf}
J.-P. Serre.
\newblock {\em Corps Locaux}.
\newblock Hermann, 1968.

\bibitem{V:book}
Helmut V{\"o}lklein.
\newblock {\em Groups as {G}alois groups}, volume~53 of {\em Cambridge Studies
  in Advanced Mathematics}.
\newblock Cambridge University Press, Cambridge, 1996.
\newblock An introduction.

\bibitem{Zhu:expsums}
H.~J. Zhu.
\newblock {$L$}-functions of exponential sums over one-dimensional affinoids:
  {N}ewton over {H}odge.
\newblock {\em Int. Math. Res. Not.}, (30):1529--1550, 2004.

\bibitem{Z:noextra}
H.~J. Zhu.
\newblock Hyperelliptic curves over {${\mathbb F}_2$} of every 2-rank without
  extra automorphisms.
\newblock {\em Proc. Amer. Math. Soc.}, 134(2):323--331 (electronic), 2006.

\end{thebibliography}

Rachel Pries\\
Colorado State University,
Mathematics department, Weber 101\\
Fort Collins, CO, 80523 (USA)\\
pries@math.colostate.edu\\

Hui June Zhu\\
SUNY at Buffalo,
Mathematics department\\
Buffalo, NY, 14260 (USA)\\
hjzhu@math.buffalo.edu
\end{document}